\documentclass[12pt,a4paper]{article}
\usepackage[T1]{fontenc}
\usepackage{lmodern}
\usepackage[margin=28mm]{geometry}
\usepackage{amsmath,amssymb,amsthm}
\usepackage{tikz,enumitem}
\usetikzlibrary{shapes.geometric}
\usepackage{float}
\usepackage[colorlinks=true,citecolor=black,linkcolor=black,urlcolor=blue]{hyperref}

\floatstyle{ruled}
\newfloat{algorithm}{htbp}{loa}
\floatname{algorithm}{Algorithm}

\newtheorem{theorem}{Theorem}
\newtheorem{lemma}[theorem]{Lemma}
\theoremstyle{definition}
\newtheorem{definition}[theorem]{Definition}

\newcommand{\PP}{\mathbb P}
\newcommand{\EE}{\mathbb E}
\newcommand{\HH}{\mathcal H}
\newcommand{\doi}[1]{\href{https://doi.org/#1}{\nolinkurl{doi:#1}}}

\title{An asymptotic version of the Burris--Schelp\\
conjecture for regular graphs of fixed degree}
\hypersetup{pdftitle={An asymptotic version of the Burris--Schelp conjecture for regular graphs of fixed degree}}
\author{Sebastian Czerwi\'nski\\
\small Institute of Mathematics, University of Zielona G\'ora, Poland\\
\small \href{mailto:s.czerwinski@im.uz.zgora.pl}{\nolinkurl{s.czerwinski@im.uz.zgora.pl}}}
\date{}

\begin{document}
\maketitle
\begin{abstract}
A proper edge-coloring is vertex-distinguishing if distinct vertices
see different sets of incident colors. For every fixed integer $d\ge2$
and every $\varepsilon>0$, we prove that every sufficiently large
$d$-regular simple graph on $n$ vertices admits a vertex-distinguishing
edge-coloring with at most $(1+\varepsilon)(d!n)^{1/d}$ colors.
The proof combines the R\"odl nibble with an entropy compression argument.
This work was prepared with the assistance of generative AI.
\end{abstract}

\section*{Declaration on the use of generative AI}
OpenAI's ChatGPT assisted with the development and presentation of the
proof. During the preparation of the manuscript, it helped identify gaps
in earlier arguments, propose ways to address them, and work out
details of the probabilistic construction, the auxiliary hypergraph, and
the entropy-compression argument. It also assisted with drafting and
revising the text, checking notation and calculations, and preparing the
\LaTeX{} source and figures. The author takes responsibility for the
mathematical claims, proofs, and final text.

\section{Introduction}
All graphs considered here are finite and simple. A graph is called
\emph{nice} if it has at most one isolated vertex and no isolated edges.
Given a (possibly partial) edge-coloring $c$, we say that a vertex $v$
\emph{sees} a color $\gamma$ if some colored edge incident with $v$
has color $\gamma$. Thus the set of colors seen by $v$ is
\[
 S_c(v)=\{c(e): e\text{ is a colored edge incident with }v\}.
\]
\begin{definition}
A proper edge-coloring of a nice graph $H$ is \emph{vertex-distinguishing}
if distinct vertices see distinct sets of colors on their incident edges.
\end{definition}
\begin{definition}
The least number of colors in a vertex-distinguishing proper edge-coloring
of $H$ is denoted by $r(H)$ and called its \emph{vertex-distinguishing index}.
\end{definition}
This parameter was studied by Burris and Schelp~\cite{burris}.
Vertex-distinguishing proper edge-colorings have also been called
\emph{strong colorings}. We use the explicit name to distinguish this
terminology from strong edge-colorings whose color classes are induced
matchings. Bazgan, Harkat-Benhamdine, Li, and Wo\'zniak~\cite{bazgan}
proved a conjecture of Burris and Schelp by showing that every nice
graph satisfies
\[
 r(H)\le |V(H)|+1.
\]
For a graph $H$, let $n_j$ be its number of vertices of degree $j$, and put
\[
 k(H)=\min\left\{k\in\mathbb N:
       \binom{k}{j}\ge n_j\text{ for every }0\le j\le\Delta(H)\right\}.
\]
Then $k(H)\le r(H)$. The stronger conjecture of Burris and Schelp asserts
that $r(H)\le k(H)+1$ for every nice graph. In the $d$-regular case,
\[
 k(H)=\min\left\{k:\binom{k}{d}\ge n\right\},
 \qquad n=|V(H)|,
\]
and in particular $r(H)\ge k$, where $k=(d!n)^{1/d}$. Graphs of maximum degree two
were studied by Balister, Bollob\'as, and Schelp~\cite{balister2}.
Balister, Kostochka, Li, and Schelp~\cite{balanced} proved the stronger
conjecture when
\[
 \Delta(H)\ge\sqrt{2|V(H)|}+4,\qquad \delta(H)\ge5.
\]

A \emph{harmonious coloring} of a simple graph $G$ is a proper vertex
coloring in which each unordered pair of colors occurs on at most one
edge. Its minimum number of colors is denoted by $h(G)$.
This notion was investigated by Frank, Harary, and Plantholt~\cite{frank}
and by Hopcroft and Krishnamoorthy~\cite{hopcroft}. If $G$ has $m$ edges,
then
\[
 \binom{h(G)}2\ge m,\qquad h(G)\ge\sqrt{2m}.
\]
Edwards~\cite{edwards} proved in 1997 that, for fixed maximum degree
and every $\varepsilon>0$, all such graphs with sufficiently many edges
satisfy
\[
 h(G)\le(1+\varepsilon)\sqrt{2m}.
\]
His proof uses the Theorem of Pippenger and Spencer and probabilistic
concentration estimates. The proof of Edwards's Theorem inspired our
proof of Theorem~\ref{main}.

For a hypergraph $\mathcal G$, a \emph{strong} or \emph{rainbow} vertex
coloring assigns distinct colors to all vertices in each edge.
It is \emph{harmonious} if distinct edges also receive distinct sets of
colors. Write $h(\mathcal G)$ for the least number of colors in such a
coloring. For a $k$-uniform hypergraph with $m$ edges,
\[
 \binom{h(\mathcal G)}k\ge m,\qquad h(\mathcal G)\ge(k!m)^{1/k}.
\]
If the rainbow requirement is dropped and only distinct edge color sets
are required, one obtains a \emph{legitimate coloring}; see
Alon and F\"uredi~\cite{alonfuredi}.

The connection with our problem is given by the dual hypergraph $H^*$.
Its vertex set is $E(H)$, and each $v\in V(H)$ defines a hyperedge
consisting of the edges incident with $v$. A proper edge-coloring of $H$
is a rainbow vertex coloring of $H^*$, and the vertex-distinguishing
condition says exactly that its hyperedges receive distinct color sets.
Consequently
\[
 r(H)=h(H^*).
\]
For a simple $d$-regular graph with $d\ge2$, the hypergraph $H^*$ is
simple, $d$-uniform, and $2$-regular. It has $n$ hyperedges.

Bosek, Czerwi\'nski, Grytczuk, and Rz\k{a}\.zewski~\cite{bcgr}
proved that a $k$-uniform hypergraph $\mathcal G$ with $m$ edges and
maximum degree $\Delta$ satisfies
\[
 h(\mathcal G)\le
 \frac{k}{k-1}\bigl(\Delta(k-1)k!m\bigr)^{1/k}+f(k,\Delta),
\]
where $f(k,\Delta)$ depends only on $k,\Delta$.
Their proof uses entropy compression and also gives a list-coloring result.
In a later single-author paper, Czerwi\'nski~\cite{czerwinski-harmonious}
gave an alternative proof of the asymptotic upper bound with the same
leading constant, using Bernshteyn's Local Cut Lemma.
D\k{e}bski, Lonc, and Rz\k{a}\.zewski~\cite{fragmentable} subsequently
obtained the asymptotically sharp bound
\[
 h(\mathcal G)\le(1+\varepsilon)(k!m)^{1/k}
\]
for fragmentable classes of bounded-degree $k$-uniform hypergraphs.
Their applications include vertex-distinguishing edge-colorings of cubic
planar graphs.

More recent results concern both general graphs and regular graphs of
growing degree. Gao, Shan, Wang, and Zhou~\cite{gaozhou} prove
\[
 r(H)\le\left\lfloor5.5k(H)+6.5\right\rfloor
\]
for every nice graph, and
\[
 r(H)\le k(H)+3
 \quad\text{if $H$ is $d$-regular and }d\ge\log_2|V(H)|\ge8.
\]
Gao, Shan, and Wang~\cite{gaoregular} determine $r(H)=d+2$ for sufficiently
large regular graphs of even order $n$ with $d$ above $n/2$ by any fixed
positive proportion. The result below concerns fixed $d$ as $n$ tends
to infinity.

\begin{theorem}\label{main}
For every fixed integer $d\ge2$ and every $\varepsilon>0$ there exists
$n_0$ such that every $d$-regular graph $H$ on $n\ge n_0$ vertices satisfies, with $k=(d!n)^{1/d}$,
\[
 r(H)\le (1+\varepsilon)k.
\]
\end{theorem}
In each stage of the proof we color an induced matching with new colors. We use Pippenger and Spencer's Theorem to choose
those colors, then complete the coloring with a separate palette.

\section{Auxiliary results}
\subsection{The Theorem of Pippenger and Spencer}
The codegree of two distinct vertices of a hypergraph is the number of
edges containing both. A proper edge-coloring assigns distinct colors
to intersecting edges.
\begin{theorem}[Pippenger--Spencer~\cite{pippenger}]\label{ps}
Fix $q\ge2$. Consider a sequence of simple $q$-uniform hypergraphs
with maximum degree $\Delta$, minimum degree $\delta$, and maximum
codegree $c$. If $\Delta\to\infty$, $\delta/\Delta\to1$, and
$c/\Delta\to0$, then $\chi'/\Delta\to1$.
\end{theorem}
Our auxiliary hypergraph contains vertices of degree one. We therefore
use the following consequence of Theorem~\ref{ps}.
\begin{lemma}\label{reg}
Fix $q\ge2$. If a simple $q$-uniform hypergraph has maximum degree at most
$\Lambda$, maximum codegree $o(\Lambda)$, and $\Lambda\to\infty$, then
its chromatic index is at most $(1+o(1))\Lambda$.
\end{lemma}
This is the asymptotic form of Corollary~16 in~\cite{katsamaktsis}.

\subsection{The bounded-differences inequality}
\begin{theorem}[McDiarmid's bounded differences inequality~\cite{mcdiarmid}]\label{bd}
Let $X_1,\dots,X_N$ be independent and let $Y=f(X_1,\dots,X_N)$.
If changing the $j$th coordinate changes $f$ by at most $c_j$, then
for every $t>0$,
\[
 \PP(|Y-\EE Y|\ge t)\le 2\exp\!\left(-\frac{2t^2}{\sum_j c_j^2}\right).
\]
\end{theorem}
This is the bounded differences form of the Azuma--Hoeffding inequality.
In the associated Doob martingale, the $j$th difference, conditional
on the preceding variables, lies in an interval of length at most $c_j$.
For a proof, see McDiarmid~\cite{mcdiarmid}.

\subsection{Chernoff bounds}
\begin{theorem}[Chernoff bounds]\label{chernoff}
Let $X$ be a sum of independent Bernoulli variables, with mean $\mu$.
For $z>0$,
\[
 \PP(X\ge\mu+z)\le\exp\!\left(-\frac{z^2}{2(\mu+z)}\right).
\]
If $u\ge\mu>0$, then
\[
 \PP(X\ge u)\le \left(\frac{\mathrm e\mu}{u}\right)^u.
\]
If $\mu=0$, then $X=0$ almost surely.
\end{theorem}
For proofs of Chernoff bounds, see Alon and Spencer~\cite[Appendix A.1]{alonspencer}
and Molloy and Reed~\cite[pp.~43--46]{molloyreed}.

\section{The proof of the main Theorem}
\subsection{Notation and introductory assumptions}
Let $(x)_r=x(x-1)\cdots(x-r+1)$, with $(x)_0=1$. Choose a large
integer $A$ depending only on $d,\varepsilon$ and perform
$A-\lceil\sqrt A\rceil$ stages. We take $\lceil\sqrt A\rceil>d$
and impose the requirement from Lemma~\ref{sampling} below.
We hold $A,d$ fixed as $n\to\infty$. All error terms in the induction
are uniform over the stages and the color sets being counted.

Put $k=(d!n)^{1/d}$, so $n=k^d/d!$, and define
\[
 \alpha_A=(A)_d^{1/d},\qquad \lambda_A=\frac{k}{\alpha_A}.
\]
For fixed $d$,
\[
 \frac{\alpha_A}{A}
 =\left(\prod_{j=0}^{d-1}\left(1-\frac{j}{A}\right)\right)^{1/d}
 \longrightarrow 1 \qquad (A\to\infty).
\]
With $A$ fixed, $\lambda_A=\Theta(k)$ as $n\to\infty$.
The palette used at each stage has size at most
$(1+o(1))\lambda_A$, as shown in \eqref{palette}. Thus the scale
of the palette is approximately $k/A$ when $A$ is large.
Since $d$ is fixed, for every
fixed exponent $\beta$ we have
\[
 n^{\beta/d}=\frac{k^\beta}{(d!)^{\beta/d}},\qquad
 O(n^{\beta/d})=O(k^\beta),\qquad o(n^{\beta/d})=o(k^\beta).
\]
We henceforth express the asymptotic estimates in terms of $k$;
the constants may depend on $d,A$. The probabilistic part uses the semi-random
method known as the \emph{R\"odl nibble}, originating in
R\"odl's packing Theorem~\cite{rodl}. It underlies the Theorem of Pippenger
and Spencer used to color the auxiliary hypergraph. For a survey of nibble
methods in graph and hypergraph colouring, see Kang, Kelly, K\"uhn,
Methuku, and Osthus~\cite{nibblesurvey}. See also Alon and
Spencer~\cite{alonspencer} and Molloy and Reed~\cite{molloyreed}.

After stage $i$, we partition the edges into the induced matchings
$E_1,\ldots,E_i$ colored in their respective stages and the residual
edge set $E(R_i)$, where $R_i$ is the residual graph. Thus $R_0=H$ and
\[
 E(H)=E_1\mathbin{\dot\cup}\cdots\mathbin{\dot\cup}E_i
       \mathbin{\dot\cup}E(R_i),
 \qquad
 E(R_i)=E(H)\setminus\bigcup_{j=1}^{i}E_j.
\]
At stage $i+1$ the selected matching $L$ becomes $E_{i+1}$, so
\[
 E_{i+1}=L,\qquad E(R_{i+1})=E(R_i)\setminus E_{i+1}.
\]
Each $E_j$ receives its own palette of new colors; the final residual
edges are colored in the $\varepsilon$-part.

After stage $i$ let $R_i$ consist of the uncolored edges and let $S_i(v)$
be the set of colors seen by $v$. Put
\[
 F_i(S)=|\{v:S_i(v)=S\}|.
\]
For distinct $(d-1)$-sets $P,Q$, let $J_i(P,Q)$ count the edges of $R_i$
whose endpoint sets are $P,Q$ in either order. We maintain:
\begin{enumerate}[label=(\roman*)]
\item Each stage colors an induced matching with a new palette.
\item All fully colored vertices are distinguished from one another.
\item For $|S|\le d-1$, with $r=d-|S|$,
\begin{equation}\label{F}
 F_i(S)\le
 \frac{(A-i)_r}{r!}\,\lambda_A^r
 +O(k^{r-1/4}).
\end{equation}
\item For distinct $(d-1)$-sets $P,Q$,
\begin{equation}\label{J}
 J_i(P,Q)\le 4\sqrt{A\lambda_A}.
\end{equation}
\end{enumerate}
We verify all four conditions at stage zero. No edge has yet been
colored, so $R_0=H$ and $S_0(v)=\varnothing$ for every vertex $v$.
Condition~(i) is vacuous, and condition~(ii) holds because every vertex
has $d\ge2$ uncolored incident edges, so there are no fully colored vertices.
For condition~(iii), the only nonempty class is the class indexed by
$S=\varnothing$. In this case $r=d$ and $k^d=d!n$, giving
\[
 F_0(\varnothing)=n
 =\frac{k^d}{d!}
 =\frac{(A)_d}{d!}\,\lambda_A^d.
\]
The last equality uses $\alpha_A^d=(A)_d$. Thus \eqref{F} holds
with equality in its main term and with zero error.
Finally, $d-1\ge1$, whereas both endpoint color sets of every edge
of $R_0$ are empty. Consequently, for all distinct $(d-1)$-sets $P,Q$,
\[
 J_0(P,Q)=0\le4\sqrt{A\lambda_A},
\]
which proves condition~(iv).
Distances between edges are measured in the line graph of $H$.
Since every color class is an induced matching in $H$, the sets
$S_i(u)$ and $S_i(v)$ are disjoint whenever $uv$ remains uncolored.
Indeed, edges of a common color at $u$ and $v$ would be at distance at
most two.

\subsection{The choice of the set \texorpdfstring{$L$}{L}}
\begin{lemma}\label{sampling}
For fixed $d$ and all sufficiently large integers $a$, every subgraph
$R\subseteq H$ admits a random induced matching $L\subseteq E(R)$
in $H$ with $\PP(e\in L)=1/a$ for each $e\in E(R)$. If $W$ is a fixed
nonempty set of vertices of residual degree $r$ and $Y_W$ counts
vertices of $W$ incident with $L$, then
\[
 \EE Y_W=\frac r a|W|,
 \qquad
 \PP(|Y_W-\EE Y_W|>t)\le 2\exp\!\left(-\frac{t^2}{C_d|W|}\right)
\]
for a constant $C_d$ depending only on $d$.
\end{lemma}
\begin{proof}
An edge has at most $2d(d-1)$ other edges at distance at most two.
Choose $p$ with $p(1-p)^{2d(d-1)}=1/a$, possible for all sufficiently
large $a$. For each edge $e$ of $R$ take an independent Bernoulli($p$)
variable. Take additional independent Bernoulli($p$) variables private
to $e$, so that their number plus the number of other residual edges
within distance two of $e$ is $2d(d-1)$. Select $e$ if its variable
is one and all these neighboring and private variables are zero.
This gives the required marginal and an induced matching in $H$.

The $r$ possible selections incident with a vertex are mutually
exclusive, giving the expectation. The count $Y_W$ depends on
$O_d(|W|)$ variables. A real edge variable affects only its own
selection and those within distance two, and a private variable affects
only one selection. Each change therefore affects $Y_W$ by $O_d(1)$.
Theorem~\ref{bd} gives the stated bound.
\end{proof}

\subsection{Stage \texorpdfstring{$i+1$}{i+1}}
Fix $0\le i<A-\lceil\sqrt A\rceil$ and put $a=A-i$. For an old
color set $S$, write $h=d-|S|$ and let $Y_S$ count the endpoints
of $L$ that see $S$. From Lemma~\ref{sampling} and \eqref{F},
simultaneously for all types with $h\ge1$ we can ensure
\begin{equation}\label{Y}
 Y_S=\frac h a F_i(S)+O(k^{h-1/4}).
\end{equation}
Indeed, $F_i(S)=O(k^h)$. For a nonempty class, the probability of
a deviation greater than $k^{h-1/4}$ is at most
\[
 2\exp(-\Omega(k^{h-1/2})).
\]
There are at most $n$ nonempty classes, so a union bound over them
gives \eqref{Y} simultaneously. Fix such an $L$.

We will color $L$ through an auxiliary $5$-uniform hypergraph $\HH_i$.
Each $f\in L$ is represented by one hyperedge. It contains two endpoint
label vertices (Section~\ref{sec:endpoint-labels}) and up to two control
vertices (Section~\ref{sec:control-blocks}). Private vertices
(Section~\ref{sec:auxiliary-hypergraph}) bring its size to five. A proper edge-coloring
of the auxiliary hypergraph $\HH_i$ prevents reuse of a color at any one label or within any one
block. The labels bound the sizes of the new classes counted by $F_{i+1}$,
and the blocks bound the repetitions counted by $J_{i+1}$.

\subsection{Control blocks for almost complete pairs}
\label{sec:control-blocks}
Consider $f=xv\in L$ and an uncolored edge $vu$, where $v$ has residual
degree two and $u$ has residual degree one. Put $U=S_i(v)$ and
$T=S_i(u)$. If $f$ receives the new color $c$, then the endpoint
sets of $vu$ change from $(U,T)$ to $(U\cup\{c\},T)$.
Both sets now have size $d-1$. Only $v$ can acquire a color in this
stage because $L$ is an induced matching in $H$.

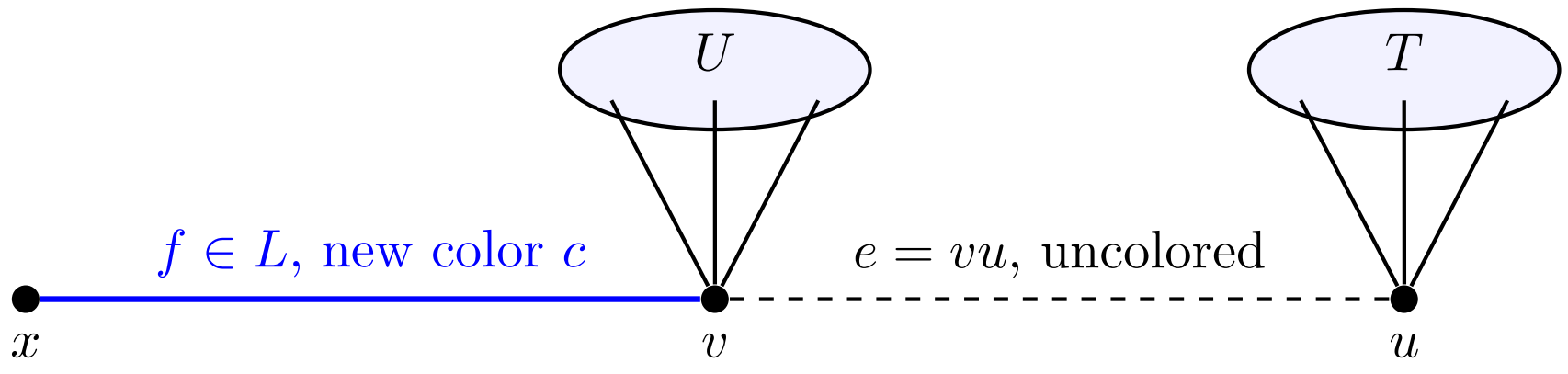
\begin{figure}[htbp]
\centering
\begin{tikzpicture}[xscale=1.55, every node/.style={font=\small}]
\node[circle,fill,inner sep=2pt,label=below:$x$] (x) at (0,0) {};
\node[circle,fill,inner sep=2pt,label=below:$v$] (v) at (3.2,0) {};
\node[circle,fill,inner sep=2pt,label=below:$u$] (u) at (6.4,0) {};
\fill[blue!5] (3.2,1.65) ellipse [x radius=0.72,y radius=0.43];
\fill[blue!5] (6.4,1.65) ellipse [x radius=0.72,y radius=0.43];
\foreach \dx in {-0.48,0,0.48} {
  \draw[thick] (v)--({3.2+\dx},1.43);
  \draw[thick] (u)--({6.4+\dx},1.43);
}
\draw[thick] (3.2,1.65) ellipse [x radius=0.72,y radius=0.43];
\draw[thick] (6.4,1.65) ellipse [x radius=0.72,y radius=0.43];
\node at (3.2,1.78) {$U$};
\node at (6.4,1.78) {$T$};
\draw[very thick,blue] (x)--node[above]{$f\in L$, new color $c$}(v);
\draw[thick,dashed] (v)--node[above]{$e=vu$, uncolored}(u);
\end{tikzpicture}
\end{figure}

Order the endpoints of each edge in $L$ as positions one and two.
For fixed $(U,T,j)$, group all incidences $(f,v)$ as above with $v$
in position $j$. Partition the group into blocks of size at most
\[
 t=\left\lceil \sqrt{A\lambda_A}\right\rceil,
\]
using the minimum possible number of blocks.
Introduce a vertex $z_B$ for each block $B$ and include it in the
hyperedge for every $f$ represented in $B$. Thus a fixed color can
occur at most once per block. The degree of $z_B$ is at most $t$,
and each $f$ belongs to at most one block per endpoint.

For fixed $U,T$, denote by $\mathcal Q_{U,T,j}$ the group of
incidences in position $j$, and put $q_j=|\mathcal Q_{U,T,j}|$.
The map sending an incidence $(f,v)$ to its corresponding vertex $u$
is injective, including across the two positions. Indeed, $u$ has only
one residual edge, which determines $vu$ and hence $v$. The vertex $v$
has exactly two residual edges, so the other one determines $f$.
Finally, the position of $v$ on $f$ is fixed. All these vertices $u$
have old set $T$. Consequently
\[
 q_1+q_2\le F_i(T).
\]
An empty group gives no blocks. Since each group is partitioned into
the minimum number of blocks of size at most $t$, the number of
blocks of type $(U,T)$ satisfies
\begin{equation}\label{blocks}
\begin{aligned}
 \left\lceil\frac{q_1}{t}\right\rceil+
 \left\lceil\frac{q_2}{t}\right\rceil
 &\le \frac{q_1+q_2}{t}+2\le \frac{F_i(T)}{t}+2\\
 &\le 2\sqrt{A\lambda_A}+2
 \le 4\sqrt{A\lambda_A}.
\end{aligned}
\end{equation}
The first inequality uses $\lceil y\rceil\le y+1$, and the second
uses the injection just described. To justify the third, apply
\eqref{F} with $|T|=d-1$, so $r=1$. Since $(A-i)/A\le1$,
\[
 F_i(T)\le
 (A-i)\lambda_A
 +O(k^{3/4})
 \le 2A\lambda_A
\]
for sufficiently large $n$. Divide by
\[
 t\ge \sqrt{A\lambda_A}.
\]
The last inequality in \eqref{blocks} holds once the displayed
quantity on the right of this lower bound is at least one.

A fixed new color $c$ occurs at most once in each block. Thus
\eqref{blocks} also bounds the number of residual edges newly acquiring
the pair $(U\cup\{c\},T)$. For a fixed resulting pair, the new color
and the side on which it occurs are uniquely determined. At most one
endpoint is hit by $L$. There are no earlier copies of that pair
because $c$ is new, while counts of pairs containing no new color
can only decrease.

\subsection{The choice of labels}
\label{sec:endpoint-labels}
For each old type $S$ occurring on an endpoint of $L$ put $r=d-1-|S|$
and introduce $B_S$ distinct labels of type $S$, where
\[
 B_S=\left\lceil\frac{(a-1)_r}{r!}\,\lambda_A^r\right\rceil.
\]
Thus $B_S=1$ when $r=0$, and $B_S=\Theta(k^r)$ otherwise.
Assign uniformly random labels to the endpoints of every edge of $L$.
Use independent assignments for distinct edges and independent
endpoint labels on an edge with different types. Equal types can occur
only when both are empty. For such an edge use a uniformly random
ordered pair of distinct labels, which is possible for large $n$ since
$d\ge2$. Each endpoint still has a uniform marginal.

Combining \eqref{F} and \eqref{Y} with $h=r+1$ gives
\begin{align*}
 Y_S &\le
 \frac{(a-1)_r}{r!}\,\lambda_A^{r+1}
       +O(k^{r+3/4})\\
 &\le B_S\lambda_A
       +O(k^{r+3/4}).
\end{align*}
The mean degree of a label is $Y_S/B_S$. Its degree is a sum
of independent Bernoulli variables indexed by the edges of $L$,
since each edge uses a particular label at most once. Theorem~\ref{chernoff}
with deviation $k^{3/4}$ gives failure probability
$\exp(-\Omega(k^{1/2}))$ for a fixed label. For $r=0$ the degree is
already deterministic. Hence simultaneously all label degrees are
at most
\[
 \lambda_A+O(k^{3/4}).
\]

For two distinct types $S,T$, the mean number of selected edges using a
fixed label of each type is at most
\[
\frac{\min(Y_S,Y_T)}{B_SB_T}=O\!\left(\frac{k}{\max(B_S,B_T)}\right).
\]
If at least one type has fewer than $d-1$ old colors this is $O(1)$.
For two labels of the empty type, each edge with two empty endpoint
sets receives that unordered pair with probability
$2/(B_\varnothing(B_\varnothing-1))$. There are at most
$Y_\varnothing/2=O(k^d)$ such edges, and
$B_\varnothing=\Theta(k^{d-1})$. Thus the mean is
$O(k^{2-d})=O(1)$.
These counts are also sums of independent indicators indexed by
edges. Theorem~\ref{chernoff} bounds their upper tails at threshold
$k^{1/2}$ by $\exp(-\Omega(k^{1/2}\log k))$.
When both types have size $d-1$, \eqref{J} gives the required bound
directly. The number of labels and pairs of labels is polynomial in
$n$, so a union bound ensures all the degree and pair bounds at once.
Fix such an assignment.

\subsection{The auxiliary hypergraph of labels and control blocks}
\label{sec:auxiliary-hypergraph}
\begin{figure}[htbp]
\centering
\begin{tikzpicture}[x=1cm,y=1cm,
 labelvertex/.style={circle,fill=black,inner sep=2pt},
 control/.style={rectangle,draw=blue!70!black,fill=blue!20,minimum size=6pt,inner sep=0pt},
 private/.style={circle,draw=black,fill=white,inner sep=2pt},
 every node/.style={font=\small}]
\draw[thick,blue!70!black,fill=blue!3,rounded corners=14pt]
 (-0.15,1.8)--(0.58,3.4)--(4.85,3.4)--(5.0,0.27)--(0.60,0.27)--cycle;
\draw[thick,teal!80!black,fill=teal!3,rounded corners=14pt]
 (0.15,1.8)--(-0.62,-0.05)--(-4.85,-0.05)--(-5.0,3.4)--(-0.62,3.4)--cycle;
\node[control,label={[label distance=3pt]left:$z_B$}] (z) at (0,1.8) {};
\node[labelvertex,label=above:$\ell_{1,x}$] (a1) at (1.2,2.15) {};
\node[labelvertex,label=above:$\ell_{1,y}$] (b1) at (3.4,2.15) {};
\node[control,label=below:$z_C$] at (1.4,1.2) {};
\node[private,label=below:$p_1$] at (3.4,1.2) {};
\node[above] at (2.75,3.46) {$h_{f_1}$};
\node[labelvertex,label=above:$\ell_{2,x}$] (a2) at (-3.5,2.15) {};
\node[labelvertex,label=above:$\ell_{2,y}$] (b2) at (-1.3,2.15) {};
\node[private,label=below:$p_2$] at (-3.3,1.05) {};
\node[private,label=below:$p'_2$] at (-1.3,1.05) {};
\node[above] at (-2.75,3.46) {$h_{f_2}$};
\draw[thick,orange!80!black,fill=orange!4,rounded corners=14pt]
 (-3.0,-0.45) rectangle (3.0,-3.35);
\node[labelvertex,label=above:$\ell_{3,x}$] (a3) at (-1.4,-1.35) {};
\node[labelvertex,label=above:$\ell_{3,y}$] (b3) at (1.4,-1.35) {};
\node[private,label=below:$p_3$] at (-1.7,-2.45) {};
\node[private,label=below:$p'_3$] at (0,-2.45) {};
\node[private,label=below:$p''_3$] at (1.7,-2.45) {};
\node[right] at (3.12,-1.85) {$h_{f_3}$};
\end{tikzpicture}
\caption{Three hyperedges of the auxiliary hypergraph $\HH_i$.
Filled circles are endpoint labels, squares are control vertices,
and open circles are private vertices. Each hyperedge $h_{f_j}$
represents an edge $f_j\in L$.}
\label{fig:auxiliary}
\end{figure}
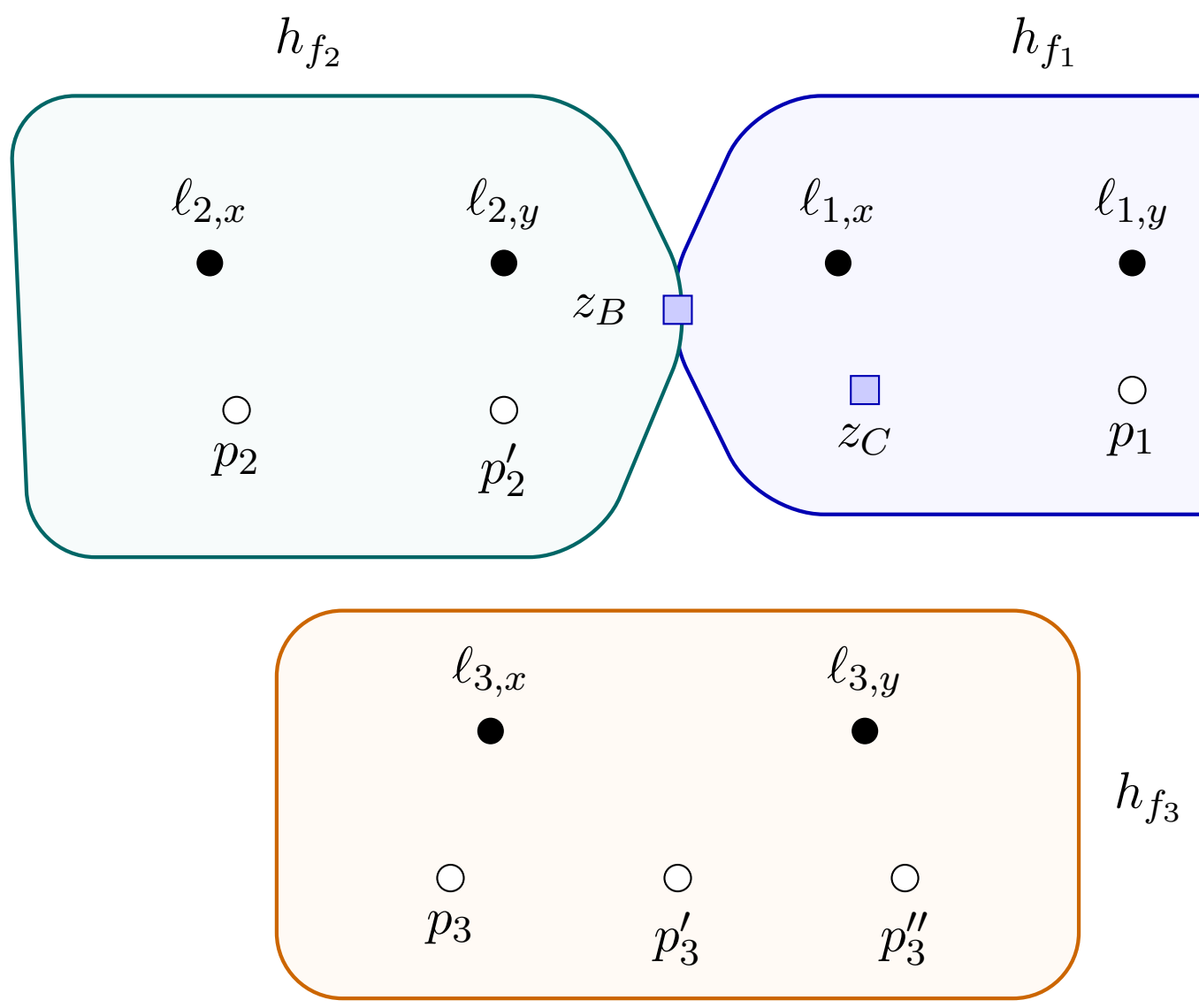

Form the auxiliary hypergraph $\HH_i$ as described above, with one hyperedge
$h_f$ for each edge $f=xy\in L$. Explicitly, $h_f$ consists of the two
endpoint labels $\ell_{f,x},\ell_{f,y}$, the zero, one, or two control
vertices of blocks containing $f$, and enough private vertices to make
its size five. Each private vertex belongs to this hyperedge only.
Figure~\ref{fig:auxiliary} shows three possible hyperedges.

The two labels of an edge are distinct. Together with at most two
control vertices they account for at most four of the five vertices
of a hyperedge. Each hyperedge therefore contains a private vertex,
so $\HH_i$ is simple. We now give the degrees
and codegrees separately for all vertex types. Write $\ell_{f,x}$ for
the label assigned to endpoint $x$ of $f$, and let $h(p)$ be the unique
hyperedge containing a private vertex $p$.

The vertex degrees are as follows.
\begin{enumerate}[label=(\roman*)]
\item For a label $\ell$ of type $S$,
\[
 \deg(\ell)=|\{(f,x):f\in L,\ x\in f,\ \ell_{f,x}=\ell\}|
 \le\lambda_A+O(k^{3/4}).
\]
No hyperedge uses the same label twice. If $|S|=d-1$, then $B_S=1$
and the degree is exactly $Y_S$.
\item For a control vertex $z_B$,
\[
 \deg(z_B)=|B|\le t
 =\left\lceil\sqrt{A\lambda_A}\right\rceil=O(k^{1/2}).
\]
There is at most one incidence of any fixed edge $f$ in $B$.
\item Every private vertex has degree exactly $1$.
\end{enumerate}
For distinct vertices $v,w$, put
$c(v,w)=|\{f\in L:v,w\in h_f\}|$. The six possible pair types are:
\begin{enumerate}[label=(\roman*)]
\item \emph{Two labels.} Their codegree is exactly the number of edges
of $L$ assigned this unordered pair of labels. For labels of distinct
types $S,T$ with $|S|=|T|=d-1$, it is at most
\[
 J_i(S,T)\le4\sqrt{A\lambda_A}.
\]
For distinct types with at least one of size at most $d-2$, and for
two labels of the empty type, the assignment chosen in Section~\ref{sec:endpoint-labels}
has codegree at most $k^{1/2}$. Two distinct labels of the same
nonempty type have codegree $0$, since the old sets on the ends of
an uncolored edge are disjoint. Thus every label pair has codegree
$O(k^{1/2})$.
\item \emph{A label and a control vertex.} Exactly
\[
 c(\ell,z_B)=|\{f\in L:\ell\in h_f,\ (f,v)\in B
                   \text{ for some }v\in f\}|
 \le |B|\le t.
\]
\item \emph{Two distinct control vertices.} Exactly
\[
\begin{aligned}
 c(z_B,z_C)&=|\{f\in L:(f,v)\in B,\ (f,w)\in C
                      \text{ for some }v,w\in f\}|
 \\ &\le\min\{|B|,|C|\}\le t.
\end{aligned}
\]
A common hyperedge uses one block at each endpoint, so here $v\ne w$.
\item \emph{A label and a private vertex.}
$c(\ell,p)=1$ if $\ell\in h(p)$, and $c(\ell,p)=0$ otherwise.
\item \emph{A control vertex and a private vertex.}
$c(z_B,p)=1$ if $z_B\in h(p)$, and $c(z_B,p)=0$ otherwise.
\item \emph{Two distinct private vertices.}
$c(p,p')=1$ if $h(p)=h(p')$, and $c(p,p')=0$ otherwise.
\end{enumerate}
Since $A,d$ are fixed and $k\to\infty$, the label degree bound
dominates the control and private degrees. Collecting these estimates gives
\begin{align*}
 \Delta(\HH_i)&\le
 \lambda_A+O(k^{3/4}),\\
 \operatorname{codeg}(\HH_i)&=O(k^{1/2}).
\end{align*}
Take a common integer degree bound
$\Lambda=\lceil\lambda_A+C k^{3/4}\rceil$ with $C$ sufficiently large.
The codegree bound is $o(\Lambda)$, so Lemma~\ref{reg} gives a proper
edge-coloring of $\HH_i$ with at most
\begin{equation}\label{palette}
 (1+o(1))\lambda_A
\end{equation}
new colors. The $o(1)$ is uniform: a sequence violating any fixed
relative error in this bound would satisfy the hypotheses of
Lemma~\ref{reg} with the same degree and codegree estimates.
Transfer the colors to the corresponding edges of $L$.

\subsection{Verification of the invariants}
Condition (i) holds by construction. To verify \eqref{F} and
the distinction of fully colored vertices in condition (ii), consider
the following cases for a color set $S$.
\begin{itemize}
\item \textbf{$S$ contains no new color.}
If $r=d-|S|\ge1$, then \eqref{Y} gives
\[
 F_{i+1}(S)=\left(1-\frac r a\right)F_i(S)+O(k^{r-1/4}).
\]
Since $(1-r/a)(a)_r=(a-1)_r$, the inductive bound on $F_i(S)$
implies \eqref{F} at stage $i+1$. If $r=0$, any vertex that sees
$S$ was already fully colored before this stage. These vertices
remain distinguished.
\item \textbf{$S$ contains a new color.}
Because $L$ is a matching, each vertex acquires at most one new
color. Thus any occurring set $S$ in this case has the form
$S=U\cup\{c\}$ for a unique new color $c$. Each label of type $U$
is used by at most one endpoint receiving $c$, so
$F_{i+1}(S)\le B_U$. If $r\ge1$, the definition of $B_U$ gives
\eqref{F} at stage $i+1$, with the rounding error absorbed in
$O(k^{r-1/4})$. If $r=0$, then $B_U=1$, so at most one newly
completed vertex has the final set $S$. Previously completed
vertices retain their old sets, which contain only old colors.
\end{itemize}

For condition (iv), fix distinct $(d-1)$-sets $P,Q$ and distinguish
two cases.
\begin{itemize}
\item \textbf{Neither $P$ nor $Q$ contains a new color.}
Every edge counted by $J_{i+1}(P,Q)$ already had these endpoint
sets before the stage. Hence $J_{i+1}(P,Q)\le J_i(P,Q)$.
\item \textbf{At least one of $P,Q$ contains a new color.}
On any edge counted by $J_{i+1}(P,Q)$, exactly one endpoint
was hit by $L$, since $L$ is induced. For a pair that occurs,
the unique new color $c$ and the set containing it are therefore
determined by $P,Q$. If $c\in P$, the selected edge at that endpoint
gives a control incidence of type $(P\setminus\{c\},Q)$. If $c\in Q$,
the type is $(Q\setminus\{c\},P)$. Each block contains at most one
selected edge colored $c$, so \eqref{blocks} gives the required bound.
\end{itemize}
This proves \eqref{J} and completes the induction.

\subsection{The \texorpdfstring{$\varepsilon$}{epsilon}-part}
The following Lemma uses the entropy-compression method of Grytczuk,
Kozik, and Micek~\cite{grytczuk}, developed for nonrepetitive sequences
and inspired by the algorithmic Lov\'asz Local Lemma of Moser and
Tardos~\cite{mosertardos}. The procedure assigns the next input color,
erases colors when a conflict occurs, and records enough information to
recover the input. Comparing the number of inputs with the number of
possible records proves that some input completes the coloring.
A related application to harmonious hypergraph colorings appears
in~\cite{bcgr}.

\begin{lemma}\label{completion}
Let $H$ be $d$-regular, with $d\ge2$, and have a proper partial edge-coloring in which fully colored
vertices are distinguished. Let $R$ consist of the uncolored edges
and let $S(v)$ be the old color set of $v$. Suppose
$S(u)\cap S(v)=\varnothing$ for every $uv\in E(R)$. Put
\[
 M_r=\max_{|S|=d-r}|\{v:S(v)=S\}|,\qquad
 M=\max(1,M_1,M_2^{1/2},\ldots,M_d^{1/d}).
\]
Then the coloring can be completed using at most $\lceil C_*M\rceil+1$
new colors, where $C_*=(d+1)(4d!+2(d-1))$.
\end{lemma}
\begin{proof}
Let $w=\lceil C_*M\rceil+1$ and fix a palette of $w$ new colors disjoint
from the old palette. For a positive integer $\ell$, let
$\mathcal C=(c_1,\ldots,c_\ell)$ be a prescribed sequence of colors
from the new palette, with repetitions allowed. Fix orders on the
vertices and residual edges. These orders determine the choice of
the next edge and the order within each local list.
For an ordered list $L$, $\operatorname{ind}_L(z)$ is the position
of $z$ in $L$, starting at~$1$. The record stores only positions
in the relevant local lists, together with case indicators and,
when needed, a permutation of local edge positions.
For an endpoint $x$ of an edge $e$, let $a(x,e)\in\{1,2\}$ be its
position among the endpoints in the vertex order. At step $j$, assign
$c_j$ to the first currently uncolored edge $e$ of $R$ and record
the outcome in $D_j$. The entries form the sequence
$\mathcal D=(D_1,D_2,\ldots)$. Stop when all edges of $R$ are colored
or the input sequence is exhausted.

The algorithm maintains the following two properties after each step.
\begin{enumerate}[label=\textup{(\roman*)}]
\item \textbf{Properness.} Any two colored edges sharing an endpoint
have different colors. An assignment at $e$ can violate this property
only if a colored edge adjacent to $e$ already has color $c_j$.
\item \textbf{Distinct sets on private edges.} For any two distinct
vertices $u,v$ with $S(u)=S(v)$ and positive degree in $R$, if all edges
in $P(v,u)\cup P(u,v)$ have received new colors, the two sets of new
colors on $P(v,u)$ and $P(u,v)$ must differ. Here
\[
 P(v,u)=E_R(v)\setminus E_R(u),
\]
where $E_R(v)$ is the set of edges of $R$ incident with $v$.
The common edge, when present, contributes equally at both vertices.
\end{enumerate}
Properness holds initially by assumption.
Condition (ii) also holds initially. Both private edge sets can be empty
only for a pair joined by an edge of $R$ with both residual
degrees one. The two old color sets have size $d-1\ge1$ and are
disjoint by hypothesis, so they cannot be equal. In the application
of Lemma~\ref{completion}, this disjointness follows from the induced matching
property of each old color class. Distinct edges of the same old color
have minimum distance at least two between their endpoints
in $H$ (equivalently, distance at least three in the line graph).
If both endpoints of a residual edge saw the same old color, the
residual edge would join endpoints of two edges of that color,
giving distance one between these endpoints.

At step $j$, write $e=vu$ for the first uncolored edge of $R$ and
color it with $c_j$. Fix an order of eligible conflicts so that each
choice below is deterministic.
\begin{itemize}
\item \textbf{Violation of (i).}
If the new color makes an improper pair, choose its common endpoint
$x\in\{u,v\}$ on $e=vu$ and the other edge $e'$. Let $L_x(e)$ be
the ordered list of residual edges at $x$ other than $e$.
The edge $e'$ belongs to this list because $c_j$ is a new color.
Write $D_j=(1,a(x,e),\operatorname{ind}_{L_x(e)}(e'))$ and erase
the color of $e$. Given $e$, there are two choices for $x$ and
at most $d-1$ positions in $L_x(e)$, so there are at most
$2(d-1)$ possible records. This case
takes priority over a violation of (ii).

\item \textbf{Violation of (ii).}
Each new equality conflict must involve \emph{one of the two
endpoints $x\in\{v,u\}$ of $e=vu$} and a vertex $y$ with
$e\in P(x,y)$. Indeed, the assignment can change
the comparison for $\{x,y\}$ only through a private edge.

Two conflicting pairs may appear at the same step, one through $v$
and one through $u$. There is at most one conflicting pair through
each endpoint in which $e$ is private to that endpoint. To see this,
suppose that $x$ has two distinct partners $y_1,y_2$. The sets
$P(x,y_1)$ and $P(x,y_2)$ together cover $E_R(x)$, since each omits
at most the edge to its respective partner. Hence all residual edges
at $x$ are colored. For each $y_i$, the private color sets agree,
and any shared residual edge is colored and contributes the same
color at both vertices. Thus both partners are fully colored and
have the same set of new colors as $x$. Their old sets also agree.
Since $e$ is private to $x$ in both pairs, neither partner is
incident with $e$. The partners therefore already violated (ii)
before this step, a contradiction.

Choose the first conflicting pair $(x,y)$ in the fixed order.
Put $b=0$ if $x,y$ are nonadjacent in $R$ and $b=1$ otherwise.
For $b=0$, let $L_x^0$ list the vertices other than $x$ that are
nonadjacent to $x$ in $R$ and have old color set $S(x)$.
For $b=1$, let $L_x^1$ list the neighbors of $x$ in $R$.
Both lists are ordered by the fixed vertex order. Properness gives
a unique color-preserving bijection $\sigma:P(x,y)\to P(y,x)$.
Write the two private edge sets in their fixed edge order as
$(f_1,\ldots,f_q)$ and $(g_1,\ldots,g_q)$, respectively, and encode
$\sigma$ by the permutation $\pi$ of $\{1,\ldots,q\}$ defined by
$\sigma(f_i)=g_{\pi(i)}$. Record
\[
 D_j=(2,a(x,e),b,\operatorname{ind}_{L_x^b}(y),\pi).
\]
Erase the colors on $P(x,y)$. Given $e$ and $q=|P(x,y)|$,
count the possible records as follows.
\begin{itemize}
\item If $x,y$ are nonadjacent in $R$, then both have residual degree
$q$. For each of the two choices of $x$, the list $L_x^0$ has
at most $M_q$ positions, and there are at most $q!$ permutations.
This gives at most $2M_q q!$ records.
\item If $x,y$ are adjacent in $R$, their old color sets are
disjoint by hypothesis. In the application, a common old color
would give two edges, one at $x$ and one at $y$, at distance at
most two in the line graph via $xy$. The induced matching
property requires distance at least three. The equality of the
two old sets thus forces both to be empty. Hence $q=d-1$.
The list $L_x^1$ has at most $d$ positions for each endpoint $x$,
and there are at most $(d-1)!$ permutations. This gives at most
$2d(d-1)!$ records. The two list indicators distinguish these
cases, so their bounds are added when both are possible.
\end{itemize}

\item \textbf{Success.} If neither condition fails, keep the color
and write $D_j=(0)$.
\end{itemize}
Every conflict created in the step requires $e$ to be colored.
In either failure case the erasure includes $e$, so it removes
all newly created conflicts, including the second equality conflict
if one occurs. Erasing colors cannot introduce a conflict.

\begin{algorithm}[htbp]
\caption{Completion by coloring, recording, and erasing.}
\label{alg:completion}
\small
\textbf{Input.} The old coloring, $R$, and
$\mathcal C=(c_1,\ldots,c_\ell)$ over $w$ new colors.\par
\textbf{Output.} A partial new coloring $\gamma$ and a record $\mathcal D$.
\begin{enumerate}[label=\arabic*.,leftmargin=8mm,itemsep=1pt,parsep=0pt,topsep=5pt]
\item Leave all edges of $R$ uncolored; set $j\leftarrow0$.
\item \textbf{while} $j<\ell$ and an edge of $R$ is uncolored \textbf{do}
\item \hspace*{3mm}$j\leftarrow j+1$.
\item \hspace*{3mm}$e\leftarrow$ the first uncolored edge of $R$.
\item \hspace*{3mm}$\gamma(e)\leftarrow c_j$.
\item \hspace*{3mm}\textbf{if} an improper pair involving $e$ occurs \textbf{then}
\item \hspace*{6mm}Choose its common endpoint $x$ and other edge $e'$.
\item \hspace*{6mm}$D_j\leftarrow(1,a(x,e),\operatorname{ind}_{L_x(e)}(e'))$.
\item \hspace*{6mm}Erase the new color of $e$.
\item \hspace*{3mm}\textbf{else if} an equality conflict occurs \textbf{then}
\item \hspace*{6mm}Choose $(x,y)$ with $e\in P(x,y)$.
\item \hspace*{6mm}Set $b$ by adjacency and encode the color bijection by $\pi$.
\item \hspace*{6mm}$D_j\leftarrow(2,a(x,e),b,\operatorname{ind}_{L_x^b}(y),\pi)$.
\item \hspace*{6mm}Erase the new colors on $P(x,y)$.
\item \hspace*{3mm}\textbf{else}
\item \hspace*{6mm}$D_j\leftarrow(0)$.
\item \hspace*{3mm}\textbf{end if}
\item \textbf{end while}
\item \textbf{return} $(\gamma,(D_1,\ldots,D_j))$.
\end{enumerate}
\end{algorithm}

The finite-input stopping condition allows the counting argument below
to consider every input sequence of length $\ell$. The old coloring is fixed
throughout. For $1\le q\le d$, summing the applicable bounds above
shows that the number of records, given $e$ and $q$, is at most
$(B_dM)^q$, where $B_d=4d!+2(d-1)$. Indeed, $M_q\le M^q$ and both
$2q!$ and $2d(d-1)!$ are at most $2d!$. The first type contributes
at most $2(d-1)$ and only when $q=1$. For $q=0$ there is one record.

\medskip
\noindent\textbf{Recovering the selected edges.}
Let $U_j$ be the set of uncolored edges of $R$ after step $j$.
We know $U_0=E(R)$. Given $U_{j-1}$, we recover $e_j$ as its
first edge in the fixed order. The entry $D_j$ then determines
the set $A_j$ of edges whose colors were erased at step $j$.
For a success $A_j=\varnothing$, while for type~1 it is $\{e_j\}$.
For type~2, the endpoint position in $D_j$ identifies $x$ on $e_j$
and $b$ selects $L_x^b$. The local position in that list identifies
$y$, so $A_j=P(x,y)$. These lists depend only on the residual graph
and the old color sets.
We therefore recover the next set by
\[
 U_j=(U_{j-1}\setminus\{e_j\})\cup A_j.
\]
This determines every $e_j$ and every $A_j$ from the records.

\medskip
\noindent\textbf{Recovering the input colors.}
Let $\gamma_j$ denote the partial coloring of $R$ after step $j$.
Starting from the final partial coloring, recover $c_j$ and
$\gamma_{j-1}$ from $\gamma_j$ and $D_j$ in reverse order.
\begin{itemize}
\item For a success, read $c_j=\gamma_j(e_j)$ and uncolor $e_j$.
\item For type~1, recover $x$ from its endpoint position and then
$e'$ from its position in $L_x(e_j)$. Read $c_j=\gamma_j(e')$,
since $e'$ retained its color. The step assigned and erased
only the color of $e_j$, so $\gamma_{j-1}=\gamma_j$.
\item For type~2, the recovered $x,y$ give the two ordered private
edge sets. The recorded permutation $\pi$ reconstructs $\sigma$
by matching their positions. Each $f\in P(x,y)$ had the color
$\gamma_j(\sigma(f))$ just before erasure. In fact, all edges
of $P(y,x)$ retained their colors because the two private edge
sets are disjoint. Thus $c_j=\gamma_j(\sigma(e_j))$.
Restore the colors of $P(x,y)\setminus\{e_j\}$ by this rule and
leave $e_j$ uncolored to obtain $\gamma_{j-1}$.
\end{itemize}
Consequently the final partial coloring and the record sequence
uniquely determine all input colors that were used.

\medskip
\noindent\textbf{Counting the records.}
Suppose no input sequence leads to completion. Then every input
of length $\ell$ uses all $\ell$ colors in the sequence, one per step.
Define its \emph{erasure-length vector} to be $(q_1,\ldots,q_\ell)$,
where $q_j=|A_j|$ counts the colors erased at step $j$, including
the color just assigned to $e_j$ whenever it is erased.
Each $q_j$ is in $\{0,\ldots,d\}$, so there are at most
$(d+1)^\ell$ such vectors. Each step adds one color to the partial
coloring before erasing $q_j$ colors. Since $R$ starts uncolored,
the number of colored edges at the end is
$\ell-\sum_{j=1}^\ell q_j\ge0$.

For a fixed vector, the prefix $D_1,\ldots,D_{j-1}$ determines $e_j$.
Hence there are at most $(B_dM)^{q_j}$ choices for $D_j$ at step $j$,
including the single choice when $q_j=0$. Multiplying these bounds
gives at most $(B_dM)^{\sum_j q_j}\le (B_dM)^\ell$ record sequences
per vector. There are at most $(w+1)^{|E(R)|}$ final partial
colorings, allowing each edge to be uncolored or to have one of the
$w$ new colors. By the reconstruction above, the $w^\ell$ inputs
give distinct pairs of a final coloring and a record sequence.
Thus
\[
 w^\ell\le (w+1)^{|E(R)|}((d+1)B_dM)^\ell.
\]
For large $\ell$ this is impossible when $w>(d+1)B_dM$.
At completion, the old palette distinguishes vertices with different
old sets, and the maintained conditions distinguish the remaining
pairs. If a pair shares a residual edge, the color of that common
edge cancels from the comparison because the coloring is proper.
\end{proof}

\subsection{The end of the proof}
After $A-\lceil\sqrt A\rceil$ stages the number of old colors is at most
\[
 (A-\lceil\sqrt A\rceil)\lambda_A+o(k)
\]
by \eqref{palette}. The hypotheses of Lemma~\ref{completion} hold by
the invariants. From \eqref{F} with $i=A-\lceil\sqrt A\rceil$,
for each $1\le r\le d$ we have
\[
 M_r\le \frac{(\lceil\sqrt A\rceil)_r}{r!}\lambda_A^r
           +O(k^{r-1/4}).
\]
Since $(\lceil\sqrt A\rceil)_r/r!\le \lceil\sqrt A\rceil^r$ and
$\lambda_A=\Theta(k)$, taking $r$th roots gives
$M_r^{1/r}\le \lceil\sqrt A\rceil\lambda_A+o(k)$. This bound also
holds for $M$ in Lemma~\ref{completion}, since $d$ is fixed.
Adding the colors from the iterative stages and from the completion,
we obtain at most
\[
 \frac{A}{\alpha_A}
 \left(1+(C_*-1)\frac{\lceil\sqrt A\rceil}{A}+o(1)\right)
 k.
\]
The coefficient apart from $o(1)$ tends to one as $A\to\infty$.
Choose $A$ first to make it less than $1+\varepsilon/2$, then
choose $n_0$ to absorb the error. This proves Theorem~\ref{main}.

\end{document}